\documentstyle[amstex,12pt,newlfont]{article}
\def \P{\Bbb{P}}
\def \Z{\Bbb{Z}}
\def \La{\Lambda}

\def \ksum{\displaystyle{\sum_{i=1}^k}}

\def \ssum{\displaystyle{\sum_{i=k+1}^{d^n}}}

\def \0{\cal{O}}
\def \s{\sigma}

\def \ed{\hspace{1cm}\Box}

\def \ni{\noindent}

\begin{document}

\title{Ample divisors on the blow up of $\P^3$ at points}

\author{Flavio Angelini}

\date{July 1996}

\maketitle

\section{Introduction}

In this note we will prove a theorem on divisors on the blow up of
$\P^3$ at points which extends a theorem of G. Xu \cite{Xu} on the
blow up of $\P^2$. The central idea of the proof works for any dimension
and therefore opens the doors to a generalization of the theorem to higher
dimension, once one overcomes certain technical difficulties that arise. 
Basically we will give a new proof of Xu's theorem
that works also for $\P^3$ and so in most of this note $n$ will be 
either 2 or 3.

\ni
{\bf Theorem.} Let $n=2$ or $3$. Let $\P^n$ be the projective space over
the field of
complex numbers. Let $p_1,\dots,p_k$ be k points in $\P^n$ in general
position and let $\pi: X \longrightarrow \P^n$ be the blow up of $\P^n$
at $p_1,\dots,p_k$ with exceptional divisors $E_1,\dots,E_k$. Let
$H=\pi^* \0_{\P^n}(1)$. Then, if $d\geq d_0(n)$, the divisor
$L=dH- \ksum E_i$ is ample if and only if $L^n >0$, i.e. $d^n>k$,
where $d_0(2)=3$, $d_0(3)=5$.

\ni
{\bf Remark 1.1.} In the case of $\P^2$ we obtain the same bound as in Xu's,
which is sharp. For $\P^3$ we believe, following a conjecture, that
the theorem holds for $d\geq3$.

Xu's proof is based on an estimate for the self-intersection of moving
singular curves in $\P^2$. The same estimate was obtained also
by Ein and Lazarsfeld in the context of Seshadri constants on smooth
surfaces \cite{EL} and used by K\"uchle to prove the above theorem in
the case of a smooth surface \cite{Ku}.

The basic idea of the proof comes from an example of R. Miranda regarding
ample divisors on a smooth surface with arbitrarily small Seshadri constant 
\cite[\S5]{L}.
The strategy is as follows:
we will use the fact that ampleness is an open condition to reduce to
the case when the points are part of the base locus of a $(n-1)$-dimensional
general linear system of hypersurfaces of degree $d$. The corresponding
line bundle
$L$ on the blow up at the entire base locus is
nef and this will imply
that the divisor obtained by blowing down some of the exceptional loci
is ample. For this last step we will use in a determinant way the fact that the
fibres of the morphism to $\P^{n-1}$ determined by $L$ are irreducible, under
the stated conditions on $d$.
We will prove this fact in Lemma 2.1 which concerns curves in $\P^3$ which are
complete intersection. The problem is to give a lower bound for the 
codimension, in the space of such curves, of the space of reducible ones.  
We emphasize that it is this technical lemma that
gives, apart from the lower bound on $d$, that may be not sharp, but
reasonable, the restriction for the dimension of the projective space for which
the theorem holds. 
It is actually possible to prove Lemma 2.1, and therefore the theorem, 
also in the cases
$n=4$ and $5$,
in a very similar way and we will spend a couple of words about it at
the end of the proof of the
lemma.
The rest of the proof works for any dimension and therefore we will present it
as much as possible in its generality in section 3.

This note is part of my Ph. D. thesis at UCLA, and I would like to
thank Rob Lazarsfeld for his guidance and encouragement.

\section{Preliminary material and lemmas}

Let $\P^N$ be the projective space parametrizing hypersurfaces in
$\P^n$ of degree $d$. So $N=N(d,n)=\binom{n+d}{n} -1$. If $F$ is such
a hypersurface, we will denote by $\left[ F \right]$ the
corresponding point in $\P^N$.
Recall that there is an action
of $PGL(N)$ on $\P^N$. If $\s$ is an element of $PGL(N)$ and $\La$ is
a linear subspace of $\P^N$, we will denote by $\La^{\s}$ the linear
subspace obtained by letting $\s$ act on $\La$. We will say that a
property holds for a general linear subspace $\La$ of $\P^N$ if it
holds for $\La^{\s}$ for $\s$ outside a union (possibly countable)
of proper subvarieties of $PGL(N)$.
We will denote
by $\Sigma_{d}$ the locus of singular hypersurfaces which
is an irreducible
subvariety of $\P^N$ of codimension one.

We will also be dealing with curves which are
complete intersection
of hypersurfaces of same degree $d$ in $\P^3$. These are
parametrized by an open subset of
the Grassmannian $Gr(\P^1,\P^N)$
(of course for $n=2$ this is just $\P^{N(d,2)}$). 
We will denote by
$\left[l_C\right] \in Gr(\P^1,\P^N)$
the point corresponding to a curve $C$.
Also we have
$dimGr(\P^k,\P^r)=(k+1)(r-k)$.
 
Let now
$$\begin{array}{rl}
NL_{d,3}=&\lbrace \text{smooth surfaces}\phantom{.} S\subset\P^3
\phantom{.}\text{of degree}\phantom{.} d \\
&\phantom{.}: Pic(S)\phantom{.} \text{is not
generated by the hyperplane class} \rbrace.
\end{array}$$
$NL_{d,3}$ is called the
{\em Noether-Lefschetz
locus} and may be viewed as a subset of $\P^{N(d,3)}$. This locus is
pretty well understood, at least as far as we are concerned here.
$NL_{d,3}$ is a countable union of
quasi-projective algebraic varieties.
The Noether-Lefschetz theorem asserts that, for $d\geq 4$, the general
surface of degree $d$ in $\P^3$ has Picard group generated by the
hyperplane section, i.e. is not contained in $NL_{d,3}$. In other words
the theorem says that the codimension of all the irreducible components
$N_i$ of
$NL_{d,3}$
is at least one.
What we will need
is an
explicit Noether-Lefschetz theorem (see \cite{Gr} for a nice proof of it)
giving a precise bound for the codimension of any of the $N_i$.

\ni
{\bf Theorem (Green).} For $d\geq 4$, the codimension of any irreducible
component of $NL_{d,3}$ in $\P^N$ is at least $(d-3)$.

The main technical lemma we need is:

\ni
{\bf Lemma 2.1.} Let $n=2,3$. Let $\P^N$ be the projective space parametrizing
hypersurfaces
of degree $d$ in $\P^n$ and assume $d\geq d_0(n)$, with $d_0(2)=3$ and
$d_0(3)=5$. Let $\La$ be a general linear subspace
of $\P^N$ of dimension $(n-1)$ whose base locus consists of $d^n$ points.
Then, for every $(n-1)$ linearly independent elements
$\left[ F_1 \right] , \dots, \left[ F_{n-1}
\right]$ of $\La$, the intersection $F_1 \cap \dots \cap F_{n-1}$ is a
curve and is irreducible.

In the case $n=2$ this is just saying that every element of a general pencil
of curves of degree $d$ is irreducible for $d\geq3$.

\ni
{\bf Remark 2.1.} Lemma 2.1, and therefore the theorem, would be proven 
for any dimension if one found
a good bound for the codimension, in the space of curves in $\P^n$ which are
a complete intersection of hypersurfaces of degree $d$, of the space of 
reducible ones. It is easy to conjecture that this codimension should be at
least $(n-1)(d-1)$, being this the codimension of such curves which have a
line as a component. This conjecture would imply Lemma 2.1 with $d_0(n)=3$ 
and for any $n$.

\ni
{\it Proof.} The proof in the case $n=2$ is rather easy. The set $\Sigma_d$
of singular curves of degree $d$ in $\P^2$ is an irreducible subvariety
of codimension one. The reducible curves are union of closed subvarieties
all lying in $\Sigma$. It is then enough to notice that, for $d\geq3$,
there are irreducible singular curves to establish that the codimension
of the reducible ones is at least two and therefore conclude that a
general pencil does not contain reducible curves.

For $n=3$ we take $d\geq5$ and we need to prove Lemma 2.1 for
a general plane $\La\subset \P^N$ with $N=\binom{d+3}{3}-1$. A curve $C$
which is a complete intersection of two elements of $\La$ corresponds to
a line
$l_C \subset\La$, i.e. to an element
$\left[l_C\right] \in Gr(\P^1,\La) \subset Gr(\P^1,\P^N)$. First
we have the following:

\ni
{\bf Claim 2.1.} For any $l_C$ in $\La$, we can find a {\em smooth}
surface $S$ with $\left[ S\right] \in l_C$ such
that $\left[ S \right]$ is not in $NL_{d,3}$.

For this we need the following sublemma:

\ni
{\bf Lemma 2.2.} For $d\geq5$, the intersection of $\La \subset\P^N$
with the Noether-Lefschetz locus consists of at most countably many points.
This is to say that there are at most countably many smooth surfaces $S$ with
$\left[ S\right] \in \La$ and with $Pic(S)$ not generated by the hyperplane
section.

\ni
{\it Proof.} We use here Kleiman's Transversality Theorem \cite[Thm. III.10.8]
{H} for the action of $PGL(N)$ on $\P^N$. We will apply the Theorem to
a plane $\La \subset \P^N$ and the closure $\overline{N_i}$ of one irreducible
component $N_i$
of $NL_{d,3}$.
The Theorem says that $\La^{\s} \cap \overline{N_i}$ is either empty or
of dimension
$$dim\La - codim\overline{N_i}$$
for $\s$ in a non-empty open subset $V_i \subset PGL(N)$.
So, by the explicit Noether-Lefschetz theorem and for $d\geq5$, we have
$$dim(\La^{\s} \cap \overline{N_i}) \leq 2-(d-3) \leq 0$$
for $\s \in V_i$.
This means that the intersection of a general plane with a component of
$NL_{d,3}$ consists at most of a finite number of points. Therefore, for $\s$
outside a possibly countable union of closed proper subvarieties, the
intersection of $\La^{\s}$ with $NL_{d,3}$ consists of at most countably many
points. $\ed$

Also observe that the intersection of $\La$ with $\Sigma_{d}$ does not
contain any line $l_C \subset \La$.
In other words, for any curve $C$ which is a complete intersection of
elements of $\La$ there exist at most finitely many surfaces $S$ with
$\left[ S\right] \in l_C$ that are singular. 
This is because $\Sigma_{d}$ is an
irreducible subvariety of codimension one and high degree of $\P^N$, and hence,
for a general $\La$,
$\Sigma_d \cap \La$ is an irreducible curve of degree strictly greater than one.

\ni
{\it Proof of Claim 2.1.} For any $l_C$ in $\La$, by Lemma 2.2
there are at most countably many surfaces $S$ with $\left[ S\right] \in
l_C$ and $Pic(S)$ not isomorphic to $\Z$ and by the observation above there
are at most finitely many singular ones. It is then possible to find
a smooth $S$ with $Pic(S)$ isomorphic to $\Z$. $\ed$

Now we can prove Lemma 2.1 for $n=3$.

\ni
{\it Proof of Lemma 2.1.} Let $C$ be any curve which is complete intersection
of elements of $\La$. We need to prove that $C$ is irreducible.
Choose a surface $S$ as in Claim 2.1.
Pick another surface $T$ with $\left[ T \right]$ in $l_C$ so that
$C=S\cap T$.
If $C$ were reducible, say $C=\cup C_i$, then any $C_i \subset S$
would be a complete intersection $S \cap T_i$, with $T_i$ a hypersurface
in $\P^3$ of degree less than $d$, since $Pic(S)$ is generated
by the hyperplane section. But then
$$C=S \cap T=\cup(S \cap T_i)=S \cap (\cup T_i)$$
This means that $\left[ \cup T_i \right] \in l_C \subset \La$.
But $\La$ misses
reducible surfaces, so $C$ has to be irreducible. $\ed$

\ni
{\bf Remark 2.2.} As mentioned in the introduction, it is possible to prove
Lemma 2.1 in the case $n=4$ and $5$ along the same line as for $n=3$, using
a generalized explicit Noether-
Lefschetz theorem, due to S. Kim \cite{Kim}, regarding smooth surfaces 
which are complete intersection of hypersurfaces \cite{A}.
The proof does not work in higher dimension due to the fact that the Noether-
Lefschetz theorem does not give any information about singular surfaces and
it is not possible anymore to ensure the existence of a smooth surface
containing the curve.  

Another tool for the proof of the theorem is the fact that
ampleness is an open condition. We will state this very well known fact
in the form of:

\ni
{\bf Proposition 2.1.} Let $\cal{L} \longrightarrow T$ be a flat family of line
bundles
over  a flat family $\cal{X} \stackrel{f}\longrightarrow T$ of
projective varieties of dimension
$n$.
Then the set
$$\lbrace t\in T\phantom{.}:\phantom{.}
\cal{L}_t=\cal{L}_{|_{X_t}}\text{is ample on}\phantom{.}X_t \rbrace$$
is open in $T$.

(See \cite{A} for a proof of it).

\section{Proof of Theorem}

The central idea of the proof works for any dimension and we will 
present the argument as much as possible in its generality. The restriction 
on the 
dimension comes 
uniquely from
Lemma 2.1, to ensure the irreducibility of the fibers of the morphism $\mu$
below.

We are given $k$ points in general position in $\P^n$ and we need to
prove that the divisor $dH-\ksum E_i$ is ample when $d^n>k$ for $d\geq d_0(n)$.
Clearly this condition is necessary.
We consider a general linear system $\La$ of hypersurfaces of degree $d$
and dimension $(n-1)$ and we let $p_1^\prime  ,\dots , p_{d^n} ^\prime$ be
the base locus of $\La$.
Let $\pi ^\prime : X^\prime \longrightarrow \P^n$ be the blow up at $k$ of
these points $p_1 ^\prime , \dots , p_k ^\prime$,
$H^\prime ={\pi ^\prime}^* \0_{\P^n}(1)$ and $E_i ^\prime, \dots ,E_k ^\prime$
the exceptional
divisors. We will prove that $dH^\prime- \ksum E_i^\prime$  is ample and
the theorem
will follow from Proposition 2.1.

To do so we consider the blow up $\nu : Z \longrightarrow \P^n$ at the whole
base locus $p_1 ^\prime, \dots , p_{d^n} ^\prime$ and set
$H=\nu^* \0_{\P^n}(1)$,
$E_i$ for $i=1, \dots ,k$, $F_i$ for $i=k+1, \dots , d^n$ the ecxeptional
divisors of $\nu$ and $L=dH-\ksum E_i -\ssum F_i$.
We have the following diagram:
$$\begin{array}{rcl}
  Z\phantom{.} &\stackrel{\delta}\longrightarrow&X^\prime \\
      &\stackrel{\nu}\searrow\phantom{.}\stackrel{\pi ^\prime}\swarrow& \\
  &\phantom{.}\P^n&
\end{array}$$
where $\delta$ is the blow down of $F_{k+1}, \dots ,F_{d^n}$.
Now $L$ is globally generated and therefore is nef. We will show that,
for $d\geq d_0(n)$, this implies that $dH^\prime-\ksum E_i^\prime$ is ample.
Now we have a morphism $\mu :Z \longrightarrow \P^{n-1}$
whose fibers are curves which are complete intersection of elements of $\La$.
Moreover, if $C$ is a curve,
$$L\cdot C=0 \Longleftrightarrow \mu(C)\phantom{..}
\mbox{\text{is a set of points}}.$$
By Lemma 2.1, for $n=2,3$ and $d\geq d_0(n)$, by taking $\La$
sufficiently general,
we may arrange for these fibers to be all irreducible.
What we need is the following:

\begin{align}
  \tag{3.1} &(L+\ssum F_i)^m \cdot Y_m \geq0 \\
  \tag{3.2} &(L+\ssum F_i)^m \cdot Y_m =0 \Longleftrightarrow Y_m \subset F_i
  \phantom{...} \mbox{\text{for one}}\phantom{.}i=k+1, \dots ,d^n,
\end{align}

where $Y_m$ is any irreducible subvariety of dimension $m$ of $Z$, for any
$0<m<n$. 

Notice that:
$$L+\ssum F_i=dH-\ksum E_i=\delta^*(dH^\prime-\ksum E_i^\prime).$$
(3.1) and (3.2) prove, by Nakai's criterion and the projection formula, that
$dH^\prime-\ksum E_i^\prime$ is ample on $X^\prime$.

To simplify notation let $L^\prime=L+\ssum F_i$. From here on we restrict to
$n=2$ and 3 because we will use Lemma 2.1. If Lemma 2.1 were proven 
for any $n$ we would proceed by induction on $m$, using the same 
arguments.
Here we just have to check (3.1) and (3.2) for $m=1$ in
the case $n=2$ and for $m=1$ and 2 in the case $n=3$. 

Let $m=1$ so that $Y_1$ is a curve. If $Y_1 \subseteq F_i$ for one $i$, then,
being $L^\prime$ trivial on $F_i$, $L^\prime \cdot Y_1 =0$ and, if
$Y_1\not\subseteq F_i$
for all $i$, then, being $L \cdot Y_1 \geq 0$ ($L$ is nef)
 and $\ssum F_i \cdot Y_1\geq 0$,
$L^\prime \cdot Y_1 \geq 0$.
This proves (3.1) and the easy direction of (3.2). For the other
direction suppose
$L^\prime \cdot Y_1=0$ but $Y_1 \not\subseteq F_i$ for all $i$. But then, being
$\ssum F_i \cdot Y_1 \geq 0$, $L \cdot Y_1$ has to be zero and therefore
$Y_1$ is a component of a fiber of $\mu$. Since 
every fiber is irreducible,  
$Y_1$ is exactly a fiber and so it meets all the exceptional divisors of $\nu$.
Therefore $\ssum F_i \cdot Y_1$ is strictly positive which is a contradiction.

In the case $n=3$ we need also to consider
any irreducible subvariety $Y_m$ of $Z$ of dimension $m=2$. As
before, if $Y_m \subseteq F_i$ for some $i$, then $L^\prime
\cdot Y_m =0$.
If $Y_m \not\subseteq F_i$ for all $i$, then we write:
$${L^\prime}^m \cdot Y_m=L \cdot {L^\prime}^{m-1}\cdot Y_m +\ssum F_i \cdot
{L^\prime}^{m-1}
\cdot Y_m.$$
Now $\ssum F_i \cdot Y_m$ has dimension strictly less than $m$ (or is empty)
and $L \cdot Y_m$ can be represented by a cycle of dimension strictly
less than $m$ (or empty) because $L$ is globally generated. So
both $L \cdot {L^\prime}^{m-1}\cdot Y_m$ and $\ssum F_i \cdot {L^\prime}^{m-1}
\cdot Y_m$
are greater or equal than zero by the previous step and so is 
${L^\prime}^m \cdot Y_m$.
For the remaining direction of (3.2) suppose ${L^\prime}^m \cdot Y_m=0$ while
$Y_m \not\subseteq F_i$ all $i$. Then, since as before $\ssum F_i \cdot
{L^\prime}^{m-1} \cdot Y_m \geq 0$, we have $L\cdot {L^\prime}^{m-1}
\cdot Y_m=0$. By the previous step, this implies that 
$L \cdot Y_m \subseteq F_i$ for
some $i$. We claim that this is a contradiction, i.e. we claim that, if
$Y_m \not\subseteq F_i$ for any $i$, then $L \cdot Y_m \not\subseteq F_i$
for any $i$. Indeed, if $Y_m \not\subseteq F_i$ for any $i$, then $\nu_*(Y_m)
\subseteq \P^3$ has dimension 2. Consider a divisor $D \in | \La|$.
By Bezout's Theorem $dim(\nu_*(Y_m) \cdot D) \geq 1$. Also
$$dim(\nu_*(Y_m) \cdot D)=dim(\nu_*(Y_m \cdot \nu^*D))=
dim(\nu_*(Y_m \cdot L)).$$
So $dim(\nu_*(Y_m\cdot L))\geq 1$ and therefore $Y_m \cdot L\not\subseteq F_i$
for any $i$. $\ed$

\ni
D\'epartement de Math\'ematiques\\Universit\'e
de Nice - Sophia-Antipolis\\Parc Valrose\\06108 Nice, France.\\
E-mail: angelini@@math.unice.fr

\end{document}